\documentclass{amsart}
\usepackage{amssymb}
\usepackage{mathrsfs}
\usepackage[T1]{fontenc}
\usepackage{xcolor}
\newcommand{\eC}{\mathscr{C}}
\newcommand{\eD}{\mathscr{D}}
\newcommand{\eE}{\mathscr{E}}
\newcommand{\eH}{\mathscr{H}}
\newcommand{\eK}{\mathscr{K}}
\newcommand{\eL}{\mathscr{L}}
\newcommand{\eM}{\mathscr{M}}
\newcommand{\eN}{\mathscr{N}}
\newcommand{\eV}{\mathscr{V}}
\newcommand{\eW}{\mathscr{W}}
\newcommand{\fD}{\mathfrak{D}}
\newcommand{\fF}{\mathfrak{F}}
\newcommand{\fL}{\mathfrak{L}}
\newcommand{\pO}{\mathcal{O}}
\newcommand{\pB}{\mathcal{B}}
\newcommand{\bC}{\mathbb{C}}
\DeclareMathOperator{\Ker}{ker}
\DeclareMathOperator{\Alg}{Alg}
\DeclareMathOperator{\Col}{Col}
\DeclareMathOperator{\Grp}{Grp}
\DeclareMathOperator{\Grf}{Grf}
\DeclareMathOperator{\Mod}{Mod}
\DeclareMathOperator{\ind}{ind}
\newtheorem{theorem}{Theorem}[section]
\newtheorem{proposition}[theorem]{Proposition}
\newtheorem{lemma}[theorem]{Lemma}

\theoremstyle{definition}

\newtheorem{example}[theorem]{Example}
\theoremstyle{remark}
\newtheorem{remark}[theorem]{Remark}
\newtheorem*{question}{Question}
\numberwithin{equation}{section}
\begin{document}
\title[Collineations of the diamond lattice]{Collineations of the diamond lattice\\ in a Hilbert space}
\author[J. Bra\v{c}i\v{c}]{Janko Bra\v{c}i\v{c}}
\address{Faculty of Natural Sciences and Engineering, University of Ljubljana, A\v{s}ker\v{c}eva c. 12, SI-1000 Ljubljana, Slovenia}
\email{janko.bracic@ntf.uni-lj.si}
\begin{abstract}
Let $\eM$ and $\eN$ be non-trivial closed subspaces of a separable Hilbert space $\eH$ such that $\eM\cap\eN=\{0\}$ and
$\eM+\eN$ is dense in $\eH$. We characterize when there exists an invertible operator on $\eH$ which maps $\eM$ onto
$\eN$ and $\eN$ onto $\eM$. We prove that such an operator exists if and only if $\eM$ and $\eN$ have a common algebraic
complement. This happens if and only if either neither of the operators $(I-P)Q$ and $(I-Q)P$ is compact, or $P-Q$ is
compact and $\dim(\eM\cap\eN^\perp)=\dim(\eM^\perp\cap\eN)$, where $P$ and $Q$ are the orthogonal projections onto
$\eM$ and $\eN$. Moreover, the operator can be chosen to be an involution. Consequently, the group of collineations of
the diamond lattice $\{\{0\},\eM,\eN,\eH\}$ which fix $\eM$ and $\eN$ is always complemented in the group of all its
collineations.
\end{abstract}
\thanks{The author was supported by the Slovenian Research and Innovation Agency program P2-0268.}
\maketitle
%
\section{Introduction}  \label{sec1}

Throughout, $\eH$ is a separable complex Hilbert space. The algebra of all bounded linear operators on $\eH$ is denoted by 
$\pB(\eH)$. A {\em subspace} is a closed linear subspace. For a subspace $\eM$, let $\eM^\perp$ be its orthogonal complement and 
let $P_{\eM}$ be the orthogonal projection onto $\eM$. For subspaces denoted by $\eM$ and $\eN$, we write $P=P_{\eM}$ and 
$Q=P_{\eN}$. Hence, $I-P=P_{\eM^\perp}$ and $I-Q=P_{\eN^\perp}$. If $\eM\subseteq\eN$, then let 
$\eN\ominus\eM=\eN\cap\eM^\perp$. For an operator $A$ on a Hilbert space 
$\eK$, let $\Grf(A)=\{x\oplus Ax;\; x\in\eK\}\subseteq\eK\oplus\eK$ be its graph.

In this note, we answer the following question.

\begin{question}
Let $\eM$ and $\eN$ be closed subspaces of $\eH$ such that $\eM\cap\eN=\{0\}$ and $\overline{\eM+\eN}=\eH$. When
does there exist an invertible operator $S\in\pB(\eH)$ such that $S\eM=\eN$ and $S\eN=\eM$?
\end{question}

Since $S$ maps $\eM$ onto $\eN$, a necessary condition is $\dim\eM=\dim\eN$. If $\eM+\eN=\eH$, that is, $\eM$ and $\eN$ are 
{\em algebraically complementary}, then this condition is also sufficient (see Lemma~\ref{lem22} below). 
For finite-dimensional $\eM$ and $\eN$, the sum $\eM+\eN$ is finite-dimensional, and therefore 
$\eM+\eN=\eH$. Hence, the interesting case of the Question is when $\dim\eM=\dim\eN=\infty$ and $\eM+\eN$ is not closed. 
The Question is natural in the setting of collineations of subspace lattices, which we now recall.

Let $\fF$ be a non-empty family of subspaces of $\eH$. Following \cite{BK22, BK26}, an invertible operator 
$S\in\pB(\eH)$ is a {\em collineation} of $\fF$ if a subspace $\eM$ belongs to $\fF$ if and only if $S\eM$ belongs 
to $\fF$. The set $\Col(\fF)$ of all collineations of $\fF$ is a group. Let $\Alg(\fF)$ be the algebra of all 
operators which leave each subspace in $\fF$ invariant. Let $\Grp(\Alg(\fF))$ be the group of all invertible 
operators $S\in\Alg(\fF)$ such that $S^{-1}\in\Alg(\fF)$. Thus, $\Grp(\Alg(\fF))$ consists of the collineations 
which fix every subspace in $\fF$. If $\fL$ is a subspace lattice, then $\Grp(\Alg(\fL))$ is a normal subgroup of 
$\Col(\fL)$ \cite[Theorem 3.7]{BK26}. A natural question is whether $\Grp(\Alg(\fL))$ is complemented in 
$\Col(\fL)$, that is, whether there is a subgroup $\pO$ of $\Col(\fL)$ such that $\Col(\fL)=\Grp(\Alg(\fL))\pO$ and 
$\Grp(\Alg(\fL))\cap\pO=\{I\}$ (see \cite[Section 3.3]{BK26}). In this case, $\Col(\fL)$ is the semidirect product 
$\Grp(\Alg(\fL))\rtimes\pO$.

We consider the simplest non-trivial lattice, the diamond. Let $\eM$ and $\eN$ be distinct non-trivial 
subspaces of $\eH$ and let $\fD=\{\{0\},\eM,\eN,\eH\}$. Then $\fD$ is a subspace lattice isomorphic to the diamond if and only if 
$\eM\cap\eN=\{0\}$ and $\overline{\eM+\eN}=\eH$. Each collineation of $\fD$ either fixes $\eM$ and $\eN$ or interchanges them (see 
Lemma~\ref{lem20}). Accordingly, let
$$\Col^{+}(\eM,\eN)=\{S\in\pB(\eH);\; S\;\text{is invertible},\; S\eM=\eM,\; S\eN=\eN\}$$
and
$$\Col^{-}(\eM,\eN)=\{S\in\pB(\eH);\; S\;\text{is invertible},\; S\eM=\eN,\; S\eN=\eM\}.$$
Note that $\Col^{+}(\eM,\eN)=\Grp(\Alg(\fD))$. Moreover, $\Col(\fD)$ is the disjoint union of 
$\Col^{+}(\eM,\eN)$ and $\Col^{-}(\eM,\eN)$. Hence, the Question asks whether $\Col^{-}(\eM,\eN)$ is non-empty.

The Question makes sense for an arbitrary pair of subspaces. It easily reduces to the case $\eM\cap\eN=\{0\}$ and 
$\overline{\eM+\eN}=\eH$. To see this, let $\eW=\eM\cap\eN$, $\eV=\overline{\eM+\eN}$, and $\eH_0=\eV\ominus\eW$. 
If $S\in\Col^{-}(\eM,\eN)$, then $S\eW=\eW$ and $S\eV=\eV$. Hence, $S_0=P_{\eH_0}S|_{\eH_0}$ is an invertible operator on $\eH_0$. 
Indeed, it is unitarily equivalent to the operator induced by $S|_{\eV}$ on $\eV/\eW$. Moreover,
$S_0\in\Col^{-}(\eM\ominus\eW,\eN\ominus\eW)$. Conversely, if $S_0\in\pB(\eH_0)$ belongs to
$\Col^{-}(\eM\ominus\eW,\eN\ominus\eW)$, then $I_{\eW}\oplus S_0\oplus I_{\eV^\perp}\in\Col^{-}(\eM,\eN)$. Note that
$(\eM\ominus\eW)\cap(\eN\ominus\eW)=\{0\}$ and that $(\eM\ominus\eW)+(\eN\ominus\eW)$ is dense in $\eH_0$.
Therefore, we may assume that $\eM\cap\eN=\{0\}$ and $\overline{\eM+\eN}=\eH$.
Subspaces with this property are called {\em topologically complementary} \cite{LP}. We also assume that $\eM$
and $\eN$ are non-trivial. Otherwise, the problem is trivial.

The answer to the Question is the following.

\begin{theorem} \label{theo11}
Let $\eM$ and $\eN$ be non-trivial topologically complementary subspaces of a separable Hilbert space $\eH$ and let
$\fD=\{\{0\},\eM,\eN,\eH\}$. The following are equivalent.
\begin{itemize}
\item[(a)] $\Col^{-}(\eM,\eN)\ne\emptyset$, that is, $\Col(\fD)\ne\Col^{+}(\eM,\eN)$.
\item[(b)] $\Col^{-}(\eM,\eN)$ contains an involution $S$, and then $\Col(\fD)=\Col^{+}(\eM,\eN)\rtimes\{I,S\}$.
\item[(c)] The subspaces $\eM$ and $\eN$ have a common algebraic complement. This means that there exists a
subspace $\eC$ such that $\eM\cap\eC=\{0\}=\eN\cap\eC$ and $\eM+\eC=\eH=\eN+\eC$.
\item[(d)] Either neither of the operators $(I-P)Q$ and $(I-Q)P$ is compact, or the operator
$P-Q$ is compact and $\dim(\eM\cap\eN^\perp)=\dim(\eM^\perp\cap\eN)$.
\end{itemize}
In particular, $\Col^{+}(\eM,\eN)$ is a complemented subgroup of $\Col(\fD)$.
\end{theorem}

Longstaff and Panaia remarked that it is not easy to prove directly that certain pairs of subspaces have no common
algebraic complement \cite[Remark 1, p.~3024]{LP}. By Theorem~\ref{theo11}, this is equivalent to the
non-existence of a collineation of $\fD$ which interchanges the subspaces.

The implication (c)$\Rightarrow$(b) in Theorem~\ref{theo11} is elementary (see Lemma~\ref{lem23} and Proposition~\ref{prop24}).
The main point is the converse implication (a)$\Rightarrow$(c). It relies on a result of Longstaff and Panaia
\cite[Theorem 2]{LP}. It also relies on the similarity invariance of two quantities attached to a pair of
subspaces. The first one is the compactness of $(I-Q)P$, and the second one is the index of a Fredholm
pair of projections \cite{ASS}. Lauzon and Treil \cite{LT} characterized common complements of two subspaces of a
Hilbert space in terms of dimensions. For the geometry of two subspaces and two projections we refer to
\cite{BS, Hal, ZH}.

The condition $\dim(\eM\cap\eN^\perp)=\dim(\eM^\perp\cap\eN)$ means that $\eM$ and $\eN$ are in the
{\em generalized generic position} in the sense of \cite{Pap}. It is sufficient for (a), but it is not necessary (see Remark~\ref{rem43}). 
If one requires $S$ to be a self-adjoint unitary involution, then the answer to the Question is different. By \cite{LWD}, such an $S$ with 
$S\eM=\eN$ exists if and only if $\eM$ and $\eN$ are in the generalized generic position.

\section{Preliminaries}  \label{sec2}

Throughout this section, let $\eM$ and $\eN$ be non-trivial topologically complementary subspaces of $\eH$, and let
$\fD=\{\{0\},\eM,\eN,\eH\}$.

\begin{lemma} \label{lem20}
The group $\Col(\fD)$ is the disjoint union of the subgroup $\Col^{+}(\eM,\eN)$ and the set $\Col^{-}(\eM,\eN)$. If
$\Col^{-}(\eM,\eN)$ is non-empty, then $\Col^{-}(\eM,\eN)=\Col^{+}(\eM,\eN)S$ for every $S\in\Col^{-}(\eM,\eN)$. In 
this case, $\Col^{+}(\eM,\eN)$ is complemented in $\Col(\fD)$ if and only if $\Col^{-}(\eM,\eN)$ contains an 
involution $S$, and then $\Col(\fD)=\Col^{+}(\eM,\eN)\rtimes\{I,S\}$.
\end{lemma}

\begin{proof}
If $S\in\pB(\eH)$ is invertible, then $S\{0\}=\{0\}$, $S\eH=\eH$, and $S\eM$, $S\eN$ are distinct non-trivial
subspaces. Hence, $S\in\Col(\fD)$ if and only if $\{S\eM,S\eN\}=\{\eM,\eN\}$ (see also \cite[Proposition 3.4]{BK26}). 
Thus, either $S\eM=\eM$ and $S\eN=\eN$, or $S\eM=\eN$ and $S\eN=\eM$. In the first case, $S\in\Col^{+}(\eM,\eN)$.
In the second case, $S\in\Col^{-}(\eM,\eN)$. Since $\eM\ne\eN$, the two cases exclude each other. Let 
$S,T\in\Col^{-}(\eM,\eN)$. Then $TS^{-1}\in\Col^{+}(\eM,\eN)$. Hence, $\Col^{-}(\eM,\eN)=\Col^{+}(\eM,\eN)S$. It 
follows that $\Col^{+}(\eM,\eN)$ has index $2$ in $\Col(\fD)$. Let $\pO$ be a complement of $\Col^{+}(\eM,\eN)$ in 
$\Col(\fD)$. Then $\pO$ is isomorphic to $\Col(\fD)/\Col^{+}(\eM,\eN)$. Hence, $\pO=\{I,S\}$, where $S^2=I$ and 
$S\not\in\Col^{+}(\eM,\eN)$, that is, $S\in\Col^{-}(\eM,\eN)$. On the other hand, if $S\in\Col^{-}(\eM,\eN)$ is an 
involution, then $\pO=\{I,S\}$ is a subgroup of $\Col(\fD)$ and $\Col^{+}(\eM,\eN)\cap\pO=\{I\}$. Moreover, 
$\Col(\fD)=\Col^{+}(\eM,\eN)\cup\Col^{+}(\eM,\eN)S=\Col^{+}(\eM,\eN)\pO$.
\end{proof}

The following canonical form is due to Longstaff and Panaia \cite[Theorem 1]{LP}. It extends Halmos' theorem
on two subspaces in generic position \cite{Hal}.

\begin{theorem}[Longstaff and Panaia] \label{theo21}
Let $\eM$ and $\eN$ be non-trivial topologically complementary subspaces of $\eH$. There exist Hilbert spaces
$\eK\ne\{0\}$ and $\eL$, with $\eL$ possibly $\{0\}$, and a positive injective contraction $H\in\pB(\eK)$ with the
following property. Possibly after interchanging $\eM$ and $\eN$, the pair $(\eM,\eN)$ is unitarily equivalent to
the pair
\begin{equation} \label{eq21}
\bigl(\Grf(-H)\oplus\eL,\;\Grf(H)\oplus\{0\}\bigr)
\end{equation}
of subspaces of $\eK\oplus\eK\oplus\eL$.\qed
\end{theorem}

If $\eM$ and $\eN$ are in the generalized generic position, then one can take $\eL=\{0\}$
\cite[Remark 1, p.~3022]{LP}. Pairs $(\eM_1,\eN_1)$ and $(\eM_2,\eN_2)$ of subspaces of Hilbert spaces $\eH_1$ and
$\eH_2$ are {\em similar} if there is an invertible operator $W\colon\eH_1\to\eH_2$ such that $W\!\!\eM_1=\eM_2$ and
$W\!\!\eN_1=\eN_2$. In this case, $W\Col^{-}(\eM_1,\eN_1)W^{-1}=\Col^{-}(\eM_2,\eN_2)$. Moreover, $S\mapsto WSW^{-1}$ maps
involutions to involutions. Note also that $\Col^{-}(\eM,\eN)=\Col^{-}(\eN,\eM)$.

\begin{lemma} \label{lem22}
Let $\eM$ and $\eN$ be non-trivial subspaces of $\eH$ such that $\eM\cap\eN=\{0\}$ and $\eM+\eN=\eH$. Then
$\Col^{-}(\eM,\eN)\ne\emptyset$ if and only if $\dim\eM=\dim\eN$. In this case, $\Col^{-}(\eM,\eN)$ contains an involution, and 
$\eM$ and $\eN$ have a common algebraic complement.
\end{lemma}

\begin{proof}
The lemma is a special case of \cite[Proposition 5.1 and Corollary 5.2]{BK26}. For the convenience of the 
reader, we include a short proof. If $S\in\Col^{-}(\eM,\eN)$, then $S|_{\eM}\colon\eM\to\eN$ is invertible. Hence, $\dim\eM=\dim\eN$. 
Conversely, let $J\colon\eM\to\eN$ be a unitary operator. The (oblique) projections of $\eH$ onto $\eM$ along $\eN$ and onto
$\eN$ along $\eM$ are bounded. Hence, $S(m+n)=J^{-1}n+Jm$ ($m\in\eM$, $n\in\eN$) defines an involution
$S\in\Col^{-}(\eM,\eN)$. The operator $m\mapsto m+Jm$ is bounded below. Therefore, the subspace
$\eC=\{m+Jm;\; m\in\eM\}$ is closed. It is easily seen that $\eM\cap\eC=\{0\}=\eN\cap\eC$ and
$\eM+\eC=\eH=\eN+\eC$.
\end{proof}

The following is part of \cite[Proposition 1]{LP}. We state it without a proof.

\begin{lemma} \label{lem23}
Let $\eM$ and $\eN$ be non-trivial topologically complementary subspaces of $\eH$ with a common algebraic complement.
Then there exist a Hilbert space $\eD$ and operators $X,Y\in\pB(\eD)$ such that the pair $(\eM,\eN)$ is similar to
the pair $\bigl(\Grf(X),\Grf(Y)\bigr)$ of subspaces of $\eD\oplus\eD$.\qed
\end{lemma}

\begin{proposition} \label{prop24}
If a pair $(\eM,\eN)$ is similar to a pair $\bigl(\Grf(X),\Grf(Y)\bigr)$ as in Lemma~\ref{lem23}, then $\Col^{-}(\eM,\eN)$
contains an involution.
\end{proposition}

\begin{proof}
It suffices to find an involution in $\Col^{-}\bigl(\Grf(X),\Grf(Y)\bigr)$. Let 
$T=\left(\begin{smallmatrix} I & 0\\ X+Y & -I\end{smallmatrix}\right)$ on $\eD\oplus\eD$. Clearly, $T^2=I$. 
Since $T(x\oplus Xx)=x\oplus Yx$ and $T(x\oplus Yx)=x\oplus Xx$ for every $x\in\eD$, we have $T\Grf(X)\subseteq\Grf(Y)$ and 
$T\Grf(Y)\subseteq\Grf(X)$. Since $T^2=I$, these inclusions are equalities.
\end{proof}

We will also need \cite[Theorem 2]{LP}. It reads as follows.

\begin{theorem}[Longstaff and Panaia] \label{theo25}
Let $\eK$ and $\eL$ be non-zero separable Hilbert spaces and let $H\in\pB(\eK)$ be a positive injective operator.
The subspaces $\Grf(-H)\oplus\eL$ and $\Grf(H)\oplus\{0\}$ of $\eK\oplus\eK\oplus\eL$ have a common algebraic
complement if and only if $H$ is non-compact.\qed
\end{theorem}

\section{Two similarity invariants}  \label{sec3}

Let $\eM$ and $\eN$ be subspaces of $\eH$. Recall that $P$ denotes the orthogonal projection onto $\eM$, and $Q$ denotes 
the orthogonal projection onto $\eN$. 

\begin{lemma} \label{lem31}
If $W\in\pB(\eH)$ is invertible, then the operator $P_{(W\!\!\eN)^\perp}P_{W\!\!\eM}$ is compact if and only if the operator $(I-Q)P$ 
is compact.
\end{lemma}

\begin{proof}
The map $J_{\eN}\colon \eN^\perp\to\eH/\eN$, which is given by $J_{\eN}y=y+\eN$, is an isometric isomorphism. Moreover,
for the quotient map $q_{\eN}\colon\eH\to\eH/\eN$, we have $q_{\eN}=J_{\eN}(I-Q)$. Hence, $(I-Q)P$ is compact if and only if 
$q_{\eN}|_{\eM}\colon\eM\to\eH/\eN$ is compact. The operator $\widehat{W}\colon \eH/\eN\to\eH/W\!\!\eN$, which is given by 
$\widehat{W}(x+\eN)=Wx+W\!\!\eN$, is invertible and $q_{W\!\!\eN}W=\widehat{W}q_{\eN}$. Thus, 
$q_{W\!\!\eN}|_{W\!\!\eM}=\widehat{W}\, q_{\eN}|_{\eM}\, W^{-1}|_{W\!\!\eM}$, where $W^{-1}|_{W\!\!\eM}\colon W\!\!\eM\to\eM$ is 
invertible. The assertion follows.
\end{proof}

\begin{lemma} \label{lem32}
The operator $P-Q$ is compact if and only if both $(I-Q)P$ and $(I-P)Q$ are compact.
\end{lemma}

\begin{proof}
If $P-Q$ is compact, then it follows from $(I-Q)(P-Q)P=(I-Q)P$ and $(I-P)(P-Q)Q=-(I-P)Q$ that $(I-Q)P$ and $(I-P)Q$ are compact.
The opposite implication follows from the equality $P-Q=((I-Q)P)^*-(I-P)Q$.
\end{proof}

Let $\pi$ be the quotient map onto the Calkin algebra. Since the kernel of $\pi$ is the ideal of compact operators,
$P-Q$ is compact if and only if $\pi(P)=\pi(Q)$. Since $\pi(P)$ and $\pi(Q)$ are projections, $(I-Q)P$ is compact if
and only if $\pi(P)\leq\pi(Q)$. Thus, Lemma~\ref{lem32} expresses the fact that $\pi(P)=\pi(Q)$ if and only if
$\pi(P)\leq\pi(Q)$ and $\pi(Q)\leq\pi(P)$.

We say that $(\eM,\eN)$ is a {\em Fredholm pair} if $P-Q$ is compact. This is a slightly more
restrictive notion than in \cite{ASS}. By Lemmas~\ref{lem31} and~\ref{lem32}, if $(\eM,\eN)$ is a Fredholm pair,
then $(W\!\!\eM,W\!\!\eN)$ is a Fredholm pair for every invertible $W\in\pB(\eH)$.

\begin{proposition} \label{prop33}
Let $(\eM,\eN)$ be a Fredholm pair. Then $T_{\eM,\eN}=Q|_{\eM}\colon \eM\to\eN$ is a Fredholm operator.
The subspaces $\eM\cap\eN^\perp$ and $\eM^\perp\cap\eN$ are finite-dimensional, and
$$\ind(\eM,\eN):=\ind\bigl(T_{\eM,\eN}\bigr)=\dim(\eM\cap\eN^\perp)-\dim(\eM^\perp\cap\eN).$$
Moreover, $\ind(\eN,\eM)=-\ind(\eM,\eN)$ and $\ind(W\!\!\eM,W\!\!\eN)=\ind(\eM,\eN)$ for every invertible
$W\in\pB(\eH)$.
\end{proposition}

\begin{proof}
Let $T=T_{\eM,\eN}$. It is easily seen that $T^*=P|_{\eN}=T_{\eN,\eM}$. For $x\in\eM$, we have
$T^*Tx=x-P(I-Q)Px$. Thus, $T^*T=I_{\eM}-K_1$, where $K_1=((I-Q)P)^*((I-Q)P)|_{\eM}$. The operator $K_1$ is compact by 
Lemma~\ref{lem32}.
Similarly, $TT^*=I_{\eN}-K_2$, where $K_2$ is compact. Hence, $T$ is invertible modulo compact operators from
both sides. By Atkinson's theorem \cite{At} (see also \cite[Chapter XI]{Con}), $T$ is a Fredholm operator. Since
$\Ker(T)=\eM\cap\eN^\perp$ and $\Ker(T^*)=\eN\cap\eM^\perp$, the formula for $\ind(\eM,\eN)$ holds. Moreover,
$\ind(\eN,\eM)=\ind(T^*)=-\ind(T)$.

Let $\eM'=W\!\!\eM$ and $\eN'=W\!\!\eN$. The range of the idempotent $E=W^{-1}P_{\eN'}W$ is $\eN$. Hence, $QE=E$ and $EQ=Q$. 
This gives $E-Q=E(I-Q)$. Therefore, $E|_{\eM}=T_{\eM,\eN}+E(I-Q)P|_{\eM}$, where the second summand is compact. It follows that
$\ind(E|_{\eM})=\ind(\eM,\eN)$. On the other hand, $E|_{\eM}=\bigl(W^{-1}|_{\eN'}\bigr)\,T_{\eM',\eN'}\,\bigl(W|_{\eM}\bigr)$, where 
$W|_{\eM}\colon\eM\to\eM'$ and $W^{-1}|_{\eN'}\colon\eN'\to\eN$ are invertible. Consequently, 
$\ind(\eM',\eN')=\ind(E|_{\eM})=\ind(\eM,\eN)$.
\end{proof}

The index $\ind(\eM,\eN)$ coincides with the index of the pair of projections $(P, Q)$ in the sense of \cite{ASS}. Now we compute both 
invariants for the canonical pair \eqref{eq21}. Statement (iii) of the following lemma is noted in \cite[Remark 2, p.~3022]{LP}.

\begin{lemma} \label{lem34}
Let $\eK\ne\{0\}$ and $\eL$ be Hilbert spaces and let $H\in\pB(\eK)$ be a positive injective contraction. Let
$\eM=\Grf(-H)\oplus\eL$ and $\eN=\Grf(H)\oplus\{0\}$. Then the following hold.
\begin{itemize}
\item[(i)] The operator $(I-P)Q$ is compact if and only if $H$ is compact.
\item[(ii)] The operator $(I-Q)P$ is compact if and only if $H$ is compact and $\eL$ is finite-dimensional.
\item[(iii)] The dimensions of the subspaces $\eM\cap\eN^\perp$ and $\eM^\perp\cap\eN$ are
$$\dim(\eM\cap\eN^\perp)=\dim\Ker(I-H)+\dim\eL\quad \text{and}\quad \dim(\eM^\perp\cap\eN)=\dim\Ker(I-H).$$
\end{itemize}
\end{lemma}

\begin{proof}
Since $H$ is self-adjoint, we have $\eM^\perp=\{Hy\oplus y\oplus 0;\; y\in\eK\}$ and
$\eN^\perp=\{(-Hy)\oplus y\oplus z;\; y\in\eK,\; z\in\eL\}$. Let $C=2(I+H^2)^{-1}H$. Then $C$ is compact
if and only if $H$ is compact.

(i) If $x\in\eK$ and $y=C x$, then $Hx-y=-H(x-Hy)$. Hence,
$x\oplus Hx\oplus 0=(Hy\oplus y\oplus 0)+\bigl((x-Hy)\oplus(-H(x-Hy))\oplus 0\bigr)\in\eM^\perp+\eM$.
Thus, $(I-P)|_{\eN}=V_+C G_+^{-1}$, where invertible operators $G_+\colon\eK\to\eN$ and $V_+\colon\eK\to\eM^\perp$
are given by $G_+x=x\oplus Hx\oplus 0$ and $V_+y=Hy\oplus y\oplus 0$, respectively.

(ii) Let $x\in\eK$, $z\in\eL$, and $y=-C x$. Then $-Hx-y=H(x+Hy)$. Hence,
$x\oplus (-Hx)\oplus z=\bigl((-Hy)\oplus y\oplus z\bigr)+\bigl((x+Hy)\oplus H(x+Hy)\oplus 0\bigr)\in\eN^\perp+\eN$.
Thus, $(I-Q)|_{\eM}=V_-\bigl((-C)\oplus I_{\eL}\bigr)G_-^{-1}$, where invertible operators $G_-\colon\eK\oplus\eL\to\eM$
and $V_-\colon\eK\oplus\eL\to\eN^\perp$ are given by $G_-(x\oplus z)=x\oplus (-Hx)\oplus z$ and 
$V_-(y\oplus z)=(-Hy)\oplus y\oplus z$, respectively.

(iii) A vector $x\oplus(-Hx)\oplus z\in\eM$ is in $\eN^\perp$ if and only if $x=H^2x$. Since $I+H$ is invertible,
this is equivalent to $x\in\Ker(I-H)$. Similarly, $x\oplus Hx\oplus 0\in\eN$ is in $\eM^\perp$ if and only if
$x\in\Ker(I-H)$.
\end{proof}

\section{Proof of the main result}  \label{sec4}

In the first two statements of this section, we assume, as in Lemma~\ref{lem34}, that $\eK\ne\{0\}$ and $\eL$ are separable Hilbert 
spaces and that $H\in\pB(\eK)$ is a positive injective contraction. Denote $\eM=\Grf(-H)\oplus\eL$ and $\eN=\Grf(H)\oplus\{0\}$.

\begin{proposition} \label{prop41}
If $H$ is non-compact, then $\Col^{-}(\eM,\eN)$ contains an involution.
\end{proposition}

\begin{proof}
Let $E(\cdot)$ be the spectral measure of $H$. Assume that $E([\varepsilon,1])\eK$ is finite-dimensional for every
$\varepsilon\in(0,1]$. Then $H$ is a norm limit of the finite-rank operators $HE([\varepsilon,1])$ as $\varepsilon\to 0$. 
This contradicts the assumption that $H$ is non-compact. Hence, there is $\varepsilon\in(0,1]$ such that $\eK_1=E([\varepsilon,1])\eK$ is 
infinite-dimensional.
Let $\eK_0=\eK\ominus\eK_1$, $H_0=H|_{\eK_0}$, and $H_1=H|_{\eK_1}$. Then $H_1\geq\varepsilon I$ is invertible.
We have $\eK\oplus\eK\oplus\eL=(\eK_0\oplus\eK_0)\oplus(\eK_1\oplus\eK_1\oplus\eL)$. With respect to this
decomposition, $\eM=\Grf(-H_0)\oplus\eM_1$ and $\eN=\Grf(H_0)\oplus\eN_1$, where $\eM_1=\Grf(-H_1)\oplus\eL$ and
$\eN_1=\Grf(H_1)\oplus\{0\}$. The unitary involution $I\oplus(-I)$ on $\eK_0\oplus\eK_0$ interchanges
$\Grf(-H_0)$ and $\Grf(H_0)$.

Note that $\eM_1\cap\eN_1=\{0\}$. We claim that $\eM_1+\eN_1=\eK_1\oplus\eK_1\oplus\eL$. Let 
$x\oplus y\oplus z\in\eK_1\oplus\eK_1\oplus\eL$. Let $u=\frac12(x-H_1^{-1}y)$ and $v=\frac12(x+H_1^{-1}y)$. Then
$x\oplus y\oplus z=\bigl(u\oplus(-H_1u)\oplus z\bigr)+\bigl(v\oplus H_1v\oplus 0\bigr)$. The subspaces $\eM_1$
and $\eN_1$ are infinite-dimensional and separable. Hence, by Lemma~\ref{lem22}, there is an involution
$S_1\in\Col^{-}(\eM_1,\eN_1)$. Thus, $(I\oplus(-I))\oplus S_1$ is an involution in $\Col^{-}(\eM,\eN)$.
\end{proof}

\begin{proposition} \label{prop42}
The following assertions are equivalent.
\begin{itemize}
\item[(a)] $\Col^{-}(\eM,\eN)\ne\emptyset$.
\item[(b)] $\Col^{-}(\eM,\eN)$ contains an involution.
\item[(c)] $H$ is non-compact or $\dim(\eM\cap\eN^\perp)=\dim(\eM^\perp\cap\eN)$.
\item[(d)] $H$ is non-compact or $\eL=\{0\}$.
\end{itemize}
\end{proposition}

\begin{proof}
(a)$\Rightarrow$(d). Let $S\in\Col^{-}(\eM,\eN)$. Assume that $H$ is compact and $\eL\ne\{0\}$. By Lemma~\ref{lem34}(i), the operator 
$(I-P)Q$ is compact. We apply Lemma~\ref{lem31} to the pair $(\eN,\eM)$ and $W=S$. It shows that 
$(I-Q)P=P_{(S\eM)^\perp}P_{S\eN}$ is compact as well. By Lemma~\ref{lem34}(ii), $\dim\eL<\infty$. By Lemma~\ref{lem32}, 
$(\eM,\eN)$ is a Fredholm pair. Now Proposition~\ref{prop33} gives $\ind(\eM,\eN)=\ind(S\eM,S\eN)=\ind(\eN,\eM)=-\ind(\eM,\eN)$.
Therefore, $\ind(\eM,\eN)=0$. On the other hand, $\Ker(I-H)$ is finite-dimensional, since $H$ is compact. Thus,
$\ind(\eM,\eN)=\dim\eL>0$, by Lemma~\ref{lem34}(iii). This is a contradiction.

(d)$\Rightarrow$(b). If $\eL=\{0\}$, then the unitary involution $I\oplus(-I)$ on $\eK\oplus\eK$ interchanges
$\eM$ and $\eN$. If $H$ is non-compact, then we apply Proposition~\ref{prop41}.

Since (b)$\Rightarrow$(a) is trivial, we proved that (a), (b), and (d) are equivalent.

(c)$\Leftrightarrow$(d). If $H$ is non-compact, then (c) and (d) both hold. Hence, we may assume that $H$ is compact, which implies 
that $\Ker(I-H)$ is finite-dimensional. By Lemma~\ref{lem34}(iii), we have $\dim(\eM\cap\eN^\perp)=\dim\Ker(I-H)+\dim\eL$ and
$\dim(\eM^\perp\cap\eN)=\dim\Ker(I-H)$. Since $\dim\Ker(I-H)<\infty$, these two dimensions are equal if and only
if $\dim\eL=0$, that is, $\eL=\{0\}$. Thus, (c) and (d) are equivalent.
\end{proof}

Now we are ready to prove the main theorem.

\begin{proof}[Proof of Theorem~\ref{theo11}]
If $\Col^{-}(\eM,\eN)$ contains an involution $S$, then $\Col(\fD)=\Col^{+}(\eM,\eN)\rtimes\{I,S\}$, by 
Lemma~\ref{lem20}. All four conditions are invariant under unitary 
equivalence and under interchanging $\eM$ and $\eN$. Hence, by
Theorem~\ref{theo21}, we may assume that $\eM=\Grf(-H)\oplus\eL$ and $\eN=\Grf(H)\oplus\{0\}$, where $\eK$ and
$\eL$ are separable.

(a)$\Leftrightarrow$(b) follows from Proposition~\ref{prop42}, and (c)$\Rightarrow$(b) follows from Lemma~\ref{lem23} and 
Proposition~\ref{prop24}.

(a)$\Rightarrow$(c). By Proposition~\ref{prop42}, $\eL=\{0\}$ or $H$ is non-compact. If $\eL=\{0\}$, then
$\{0\}\oplus\eK$ is a common algebraic complement of $\Grf(-H)$ and $\Grf(H)$. If $\eL\ne\{0\}$ and $H$ is
non-compact, then $\eM$ and $\eN$ have a common algebraic complement by Theorem~\ref{theo25}.

(a)$\Leftrightarrow$(d). By Proposition~\ref{prop42}, (a) holds if and only if $H$ is non-compact or $\eL=\{0\}$. We show that this is 
equivalent to (d). If $H$ is non-compact, then neither $(I-P)Q$ nor $(I-Q)P$ is compact, by Lemma~\ref{lem34}. Let $H$ be compact and 
$\eL=\{0\}$. Then $P-Q$ is compact, by Lemmas~\ref{lem32} and~\ref{lem34}. The dimensions in (d) are equal, by
Lemma~\ref{lem34}(iii). Conversely, assume (d). If neither $(I-P)Q$ nor $(I-Q)P$ is compact, then $H$ is non-compact. Let 
$P-Q$ be compact. Then $H$ is compact and $\dim\eL<\infty$, by Lemmas~\ref{lem32} and~\ref{lem34}. Hence, $\Ker(I-H)$ is 
finite-dimensional. Now the equality of dimensions gives $\eL=\{0\}$, by Lemma~\ref{lem34}(iii).

If $\Col^{-}(\eM,\eN)=\emptyset$, then $\Col(\fD)=\Col^{+}(\eM,\eN)$, and $\{I\}$ is a complement of 
$\Col^{+}(\eM,\eN)$. On the other hand, if $\Col^{-}(\eM,\eN)\ne\emptyset$, then $\Col^{-}(\eM,\eN)$ contains an 
involution $S$, by (b). Hence, $\{I,S\}$ is a complement of $\Col^{+}(\eM,\eN)$ in $\Col(\fD)$.
\end{proof}

\begin{remark} \label{rem43}
(i) The condition $\dim(\eM\cap\eN^\perp)=\dim(\eM^\perp\cap\eN)$ is not necessary for (a). For instance, let
$\eK=\ell^2$, $\eL=\bC$, $H=\frac12 I$, $\eM=\Grf(-H)\oplus\eL$, and $\eN=\Grf(H)\oplus\{0\}$. Then
$\dim(\eM\cap\eN^\perp)=1$ and $\eM^\perp\cap\eN=\{0\}$. However, $H$ is non-compact. Hence,
$\Col^{-}(\eM,\eN)\ne\emptyset$. However, by \cite{LWD}, there is no self-adjoint unitary involution $S$ with
$S\eM=\eN$. On the other hand, let $H$ be compact. Then $\Col^{-}(\eM,\eN)\ne\emptyset$ if and only
if $\eM$ and $\eN$ are in the generalized generic position.

(ii) Let $H$ be compact and let $\eL$ be infinite-dimensional. Then exactly one of the operators
$(I-P)Q$ and $(I-Q)P$ is compact. Hence, $\Col^{-}(\eM,\eN)=\emptyset$ follows from
Lemma~\ref{lem31} alone. Together with Lemma~\ref{lem23} and Proposition~\ref{prop24}, this gives a short proof
that such $\eM$ and $\eN$ have no common algebraic complement (cf.\ \cite[Remark 1, p.~3024]{LP}).\qed
\end{remark}

The following example is due to Professor Vladim\'{i}r M\"{u}ller. The author is thankful for permission to include it in the paper. The 
example shows that in some cases, $\Col^{-}(\eM,\eN)=\emptyset$ can be proved by an elementary argument.

\begin{example} \label{ex44}
Let $\eH$ be a Hilbert space with an orthonormal basis $\{e_j,f_j,g_j;\; j\in\mathbb{N}\}$. Let
$\eM=\bigvee\{e_j,\, f_j+2^{-j}g_j;\; j\in\mathbb{N}\}$ and $\eN=\bigvee\{f_j;\; j\in\mathbb{N}\}$.
It is easily seen that $\eM\cap\eN=\{0\}$ and $\overline{\eM+\eN}=\eH$. Moreover, the subspaces $\eM$, $\eN$,
$\eM^\perp$, and $\eN^\perp$ are infinite-dimensional. We show that $\Col^{-}(\eM,\eN)=\emptyset$.

Assume that there exists $S\in\Col^{-}(\eM,\eN)$. Let $\eE=\bigvee\{e_j;\; j\in\mathbb{N}\}$ and, for $n\in\mathbb{N}$,
let $\eN_n=\bigvee\{f_j;\; j\geq n\}$. Let $V\colon\eN\to\eM$ be the bounded operator given by
$Vf_j=f_j+2^{-j}g_j$ ($j\in\mathbb{N}$). Since $\eE\subseteq\eM$, we have $S^{-1}\eE\subseteq S^{-1}\eM=\eN$. The
subspace $S^{-1}\eE$ is infinite-dimensional and $\eN_n$ has finite codimension in $\eN$. Hence, there exists a unit
vector $x\in\eN_n\cap S^{-1}\eE$. Then $\|Vx-x\|\leq 2^{-n}$. Moreover, $Sx\in \eE$ and $SVx\in S\eM=\eN$. Since
$\eE\perp\eN$, we have $2^{-n}\|S\|\geq\|S(Vx-x)\|\geq\|Sx\|\geq\|S^{-1}\|^{-1}$.
This holds for every $n\in\mathbb{N}$, which is not possible. Thus, $\Col^{-}(\eM,\eN)=\emptyset$. By
Theorem~\ref{theo11}, the subspaces $\eM$ and $\eN$ have no common algebraic complement.

This also follows from Remark~\ref{rem43}(ii). Indeed, since $\eE\subseteq\eM\cap\eN^\perp$, the operator $(I-Q)P$
is not compact. On the other hand, $\eM^\perp=\bigvee\{g_j-2^{-j}f_j;\; j\in\mathbb{N}\}$ and
$(I-P)f_j=-\frac{2^{-j}}{1+4^{-j}}(g_j-2^{-j}f_j)$. Therefore, $(I-P)Q$ is compact. The
argument above is elementary, but it uses the fact that $\eM\cap\eN^\perp$ is infinite-dimensional. If
$\eM\cap\eN^\perp$ is finite-dimensional, then the index of Proposition~\ref{prop33} seems to be needed.\qed
\end{example}

\section{Operators mapping $\eM$ into $\eN$ and $\eN$ into $\eM$}  \label{sec5}

Throughout this section, let $\eM$ and $\eN$ be non-trivial topologically complementary subspaces of $\eH$ and let
$\fD=\{\{0\},\eM,\eN,\eH\}$. Recall that $\Alg(\fD)=\{A\in\pB(\eH);\; A\eM\subseteq\eM,\; A\eN\subseteq\eN\}$. Let
$$\Mod(\fD)=\{T\in\pB(\eH);\; T\eM\subseteq\eN,\; T\eN\subseteq\eM\}.$$
The definition of $\Mod(\fD)$ is symmetric in $\eM$ and $\eN$. Hence, $\Mod(\fD)$ depends only on the lattice
$\fD$, as $\Alg(\fD)$ does. Thus, the notation $\Mod(\fD)$ is legitimate. The sets $\Col^{\pm}(\eM,\eN)$ also
depend only on $\fD$. We keep the notation $\Col^{\pm}(\eM,\eN)$ from the previous sections.

It is clear that $\Alg(\fD)$ is an algebra and that $\Mod(\fD)$ is a subspace of $\pB(\eH)$. Both are closed in the
weak operator topology. It is easily seen that
$$\Alg(\fD)\Mod(\fD)\subseteq\Mod(\fD),\quad \Mod(\fD)\Alg(\fD)\subseteq\Mod(\fD),\quad
\Mod(\fD)^2\subseteq\Alg(\fD).$$
If $A\in\Alg(\fD)\cap\Mod(\fD)$, then $A\eM\subseteq\eM\cap\eN=\{0\}$ and $A\eN\subseteq\eN\cap\eM=\{0\}$. Since
$\eM+\eN$ is dense in $\eH$, we have $A=0$. Thus, $\Alg(\fD)\oplus\Mod(\fD)$ is a $\mathbb{Z}_2$-graded algebra.
The group $\Col^{+}(\eM,\eN)$ consists of the invertible operators $A\in\Alg(\fD)$ such that
$A^{-1}\in\Alg(\fD)$. Similarly, $\Col^{-}(\eM,\eN)$ consists of the invertible operators $T\in\Mod(\fD)$ such that
$T^{-1}\in\Mod(\fD)$. Indeed, if $T,T^{-1}\in\Mod(\fD)$, then $\eN=TT^{-1}\eN\subseteq T\eM\subseteq\eN$ and,
similarly, $T\eN=\eM$. Example~\ref{ex53} below shows that the condition $T^{-1}\in\Mod(\fD)$ cannot be omitted.

We need the following lemma, which is folklore.

\begin{lemma} \label{lem50}
Let $\eM\subseteq\eH$ be a subspace. If $T\in\pB(\eH)$ is invertible, then $T^*\eM^\perp=(T^{-1}\eM)^\perp$.
\end{lemma}

\begin{proof}
For $x\in\eH$, we have $x\perp T^*\eM^\perp$ if and only if $Tx\perp\eM^\perp$, that is, $Tx\in\eM$. Hence,
$(T^*\eM^\perp)^\perp=T^{-1}\eM$. Since $T^*$ is invertible, $T^*\eM^\perp$ is closed, and therefore
$T^*\eM^\perp=(T^{-1}\eM)^\perp$.
\end{proof}

\begin{proposition} \label{prop51}
For an invertible operator $W\in\pB(\eH)$, the following assertions are equivalent.
\begin{itemize}
\item[(a)] $\Mod(\fD)=W\Alg(\fD)$.
\item[(b)] $\Mod(\fD)=\Alg(\fD)W$.
\item[(c)] $W\in\Col^{-}(\eM,\eN)$.
\end{itemize}
In particular, $\Mod(\fD)$ is a translate of $\Alg(\fD)$ if and only if $\eM$ and $\eN$ have a common algebraic
complement. In this case, there exists an involution $S\in\Col^{-}(\eM,\eN)$ such that $\Mod(\fD)=S\Alg(\fD)$ and
$S\Alg(\fD)S=\Alg(\fD)$.
\end{proposition}

\begin{proof}
Put $\fD^\perp=\{\{0\},\eM^\perp,\eN^\perp,\eH\}$. For every $T\in\pB(\eH)$, we have $T\eM\subseteq\eN$ if and only if
$T^*\eN^\perp\subseteq\eM^\perp$. Hence, $\Mod(\fD)^*=\Mod(\fD^\perp)$ and $\Alg(\fD)^*=\Alg(\fD^\perp)$. Moreover, if
$W$ is invertible, then $(W\!\!\eM)^\perp=(W^*)^{-1}\eM^\perp$, by Lemma~\ref{lem50} applied to $T=W^{-1}$. 
Thus, $W\!\!\eM=\eN$ if and only if 
$W^*\eN^\perp=\eM^\perp$, and therefore $W\in\Col^{-}(\eM,\eN)$ if and only if $W^*\in\Col^{-}(\eM^\perp,\eN^\perp)$. Taking 
adjoints, we see that (b) for $W$ and $\fD$ is (a) for $W^*$ and $\fD^\perp$. The subspaces $\eM^\perp$ and $\eN^\perp$ are 
non-trivial, $\eM^\perp\cap\eN^\perp=(\eM+\eN)^\perp=\{0\}$, and $\eM^\perp+\eN^\perp$ is dense, since 
$(\eM^\perp+\eN^\perp)^\perp=\eM\cap\eN=\{0\}$. Hence, $\fD^\perp$ satisfies the same assumptions as $\fD$. Therefore, it 
suffices to prove that (a) and (c) are equivalent.

(c)$\Rightarrow$(a). Let $W\in\Col^{-}(\eM,\eN)$. If $A\in\Alg(\fD)$, then $WA\eM\subseteq W\!\!\eM=\eN$ and
$WA\eN\subseteq W\!\!\eN=\eM$. Hence, $WA\in\Mod(\fD)$. Conversely, let $T\in\Mod(\fD)$. Since $W^{-1}\eN=\eM$ and
$W^{-1}\eM=\eN$, the operator $W^{-1}T$ belongs to $\Alg(\fD)$. Hence, $T=W(W^{-1}T)\in W\Alg(\fD)$.

(a)$\Rightarrow$(c). Since $I\in\Alg(\fD)$, we have $W\in\Mod(\fD)$, that is, $W\!\!\eM\subseteq\eN$ and
$W\!\!\eN\subseteq\eM$. For $x,y\in\eH$, let $x\otimes y$ be the operator $v\mapsto\langle v,y\rangle x$. Since $\eM$ is
non-trivial and $\eM\cap\eN=\{0\}$, we have $\eN\ne\eH$, and therefore there exists a non-zero $y\in\eN^\perp$. Let $x\in\eN$.
Then $x\otimes y\in\Mod(\fD)$, and therefore $W^{-1}(x\otimes y)=(W^{-1}x)\otimes y$ maps $\eM$ into $\eM$. Since
$\eM^\perp\cap\eN^\perp=\{0\}$, there exists $m\in\eM$ with $\langle m,y\rangle\ne0$. Hence, $W^{-1}x\in\eM$. This shows that
$\eN\subseteq W\!\!\eM$. Interchanging the roles of $\eM$ and $\eN$, we obtain $\eM\subseteq W\!\!\eN$. Thus,
$W\in\Col^{-}(\eM,\eN)$.

By the equivalence of (a) and (c), $\Mod(\fD)$ is a translate of $\Alg(\fD)$ if and only if
$\Col^{-}(\eM,\eN)\ne\emptyset$. By Theorem~\ref{theo11}, this happens if and only if $\eM$ and $\eN$ have a common
algebraic complement, and in this case $\Col^{-}(\eM,\eN)$ contains an involution $S$. If $A\in\Alg(\fD)$, then
$SAS\eM=SA\eN\subseteq S\eN=\eM$ and, similarly, $SAS\eN\subseteq\eN$. Hence, $S\Alg(\fD)S\subseteq\Alg(\fD)$. Since
$S^2=I$, we also have $\Alg(\fD)=S\bigl(S\Alg(\fD)S\bigr)S\subseteq S\Alg(\fD)S$.
\end{proof}

Now we show that $\Mod(\fD)$ can contain invertible operators even if $\Col^{-}(\eM,\eN)$ is empty. 
As in Section~\ref{sec4}, let $\eK\ne\{0\}$ and $\eL$ be separable Hilbert spaces, and let $\eH=\eK\oplus\eK\oplus\eL$.
Let $H\in\pB(\eK)$ be a positive injective contraction. In $\eH$ we consider the subspaces $\eM=\Grf(-H)\oplus\eL$ and 
$\eN=\Grf(H)\oplus\{0\}$. For $A\in\pB(\eK)$ and $C\in\pB(\eL,\eK)$, let $[A\ C]$ be the operator $x\oplus z\mapsto Ax+Cz$ from 
$\eK\oplus\eL$ to $\eK$. For $D\in\pB(\eK)$ and $E\in\pB(\eK,\eL)$, let $\left[\begin{smallmatrix} D\\ E\end{smallmatrix}\right]$ be the
operator $x\mapsto Dx\oplus Ex$ from $\eK$ to $\eK\oplus\eL$.

\begin{proposition} \label{prop52}
An operator on $\eK\oplus\eK\oplus\eL$ belongs to $\Mod(\fD)$ if and only if it is of the form
\begin{small}
$$T=\left(\begin{array}{ccc} A & B & C\\ -HBH & D & HC\\ EH & E & 0\end{array}\right),$$
\end{small}
where $A,B,D\in\pB(\eK)$, $C\in\pB(\eL,\eK)$, $E\in\pB(\eK,\eL)$, and $HA=-DH$.

Let $T\in\Mod(\fD)$ be as above with $B=0$. Then $T$ is invertible if and only if the operators $[A\ C]$ and
$\left[\begin{smallmatrix} D\\ E\end{smallmatrix}\right]$ are invertible.
\end{proposition}

\begin{proof}
Let $T=\left(\begin{smallmatrix} A & B & C\\ U & D & V\\ W & E & Z\end{smallmatrix}\right)$ be an operator on
$\eK\oplus\eK\oplus\eL$. For $x\in\eK$ and $z\in\eL$, we have
\begin{align*}
T(x\oplus(-Hx)\oplus z)&=\bigl((A-BH)x+Cz\bigr)\oplus\bigl((U-DH)x+Vz\bigr)\\
&\qquad\oplus\bigl((W-EH)x+Zz\bigr),\\
T(x\oplus Hx\oplus 0)&=(A+BH)x\oplus(U+DH)x\oplus(W+EH)x.
\end{align*}
The first vector is in $\eN$ for all $x$ and $z$ if and only if $U-DH=H(A-BH)$, $V=HC$, $W=EH$, and $Z=0$. The second
vector is in $\eM$ for all $x$ if and only if $U+DH=-H(A+BH)$. These two equalities are equivalent to $U=-HBH$ and
$HA=-DH$.

Now let $B=0$. Let $F\colon\eK\oplus\eK\oplus\eL\to\eK\oplus\eK$ and
$G\colon\eK\oplus\eK\to\eK\oplus\eK\oplus\eL$ be the operators given by
$$F(x\oplus y\oplus z)=(Ax+Cz)\oplus(Hx+y)\qquad\text{and}\qquad G(x\oplus y)=x\oplus(Hx+Dy)\oplus Ey.$$
Since $HA+DH=0$, a direct computation gives $T=GF$. Assume that $[A\ C]$ and
$\left[\begin{smallmatrix} D\\ E\end{smallmatrix}\right]$ are invertible. Let $u\oplus v\in\eK\oplus\eK$ and let
$x\oplus z=[A\ C]^{-1}u$. Then $F(x\oplus(v-Hx)\oplus z)=u\oplus v$. If $F(x\oplus y\oplus z)=0$, then
$Ax+Cz=0$ and $y=-Hx$. Hence, $x=0$, $z=0$, and $y=0$. Thus, $F$ is invertible. Similarly, let
$x\oplus y\oplus z\in\eK\oplus\eK\oplus\eL$ and let
$u=\left[\begin{smallmatrix} D\\ E\end{smallmatrix}\right]^{-1}\bigl((y-Hx)\oplus z\bigr)$. Then
$G(x\oplus u)=x\oplus y\oplus z$, and $G$ is injective. Thus, $G$ is invertible. Therefore, $T=GF$ is invertible.

Conversely, assume that $T$ is invertible. The first coordinate of $T(x\oplus y\oplus z)$ is $Ax+Cz$. Hence, $[A\ C]$
is surjective. Since $HA=-DH$, we have $T(x\oplus(-Hx)\oplus z)=(Ax+Cz)\oplus H(Ax+Cz)\oplus 0$. If $Ax+Cz=0$, then
$T(x\oplus(-Hx)\oplus z)=0$, and therefore $x=0$ and $z=0$. Thus, $[A\ C]$ is injective. Since
$T(0\oplus y\oplus 0)=0\oplus Dy\oplus Ey$, the operator $\left[\begin{smallmatrix} D\\ E\end{smallmatrix}\right]$ is
injective. Let $s\in\eK$ and $t\in\eL$, and let $T(x\oplus y\oplus z)=0\oplus s\oplus t$. Then $Ax+Cz=0$. Hence, $x=0$
and $z=0$, and therefore $Dy=s$ and $Ey=t$. Thus, $\left[\begin{smallmatrix} D\\ E\end{smallmatrix}\right]$ is
surjective as well.
\end{proof}

The following example shows that $\Mod(\fD)$ can contain an invertible operator even if $\Col^{-}(\eM,\eN)$ is empty.

\begin{example} \label{ex53}
We use the notation of Proposition~\ref{prop52} with $\eL=\bC$. Let $(e_k)_{k\geq 0}$ be an orthonormal basis of $\eK$ and let 
$He_k=h_ke_k$, where $h_k=1/k!$. Then $H$ is a compact positive injective contraction. Let $R$ be the unilateral shift, $Re_k=e_{k+1}$, 
and let $U$ and $V$ be the compact weighted shifts given by $Ue_0=0$, $Ue_k=\frac{1}{k}\,e_{k-1}$ ($k\geq 1$), and 
$Ve_k=\frac{1}{k+1}\,e_{k+1}$ ($k\geq 0$). Then $X=R+U$ and $Y=R^*+V$ are Fredholm operators of index $-1$ and $1$, respectively. 
Since $h_k/h_{k-1}=1/k$, we have $YHe_k=h_ke_{k-1}+h_{k+1}e_{k+1}=HXe_k$ for every $k\geq 0$, where $e_{-1}$ is meant to be $0$. 
Hence, $YH=HX$.

Let $Xx=0$, where $x=\sum_k x_ke_k$. Comparing the coefficients, we get $x_1=0$ and $x_{m+1}=-(m+1)x_{m-1}$ for
$m\geq 1$. Hence, $|x_{2j}|=2^jj!\,|x_0|$. Since $x\in\eK$, we get $x=0$. Thus, $X$ is injective and its range has
codimension $1$. Similarly, $Yx=0$ if and only if $x_1=0$ and $x_{m+1}=-\frac{1}{m}\,x_{m-1}$ for $m\geq 1$. Hence,
$\Ker(Y)$ is spanned by $u=\sum_{j\geq 0}\frac{(-1)^j}{1\cdot 3\cdots(2j-1)}\,e_{2j}$, and $Y$ is surjective.

Let $A=X$ and $D=-Y$, let $C\in\pB(\bC,\eK)$ map $1$ to a vector which is not in the range of $X$, and let
$Ex=\langle x,e_0\rangle$. Then $HA=-DH$, and $[A\ C]$ is invertible. Since $Eu=1$, the operator
$\left[\begin{smallmatrix} D\\ E\end{smallmatrix}\right]$ is invertible as well. By Proposition~\ref{prop52},
$\Mod(\fD)$ contains an invertible operator $T$. Since $H$ is compact and $\eL\ne\{0\}$, we have $\Col^{-}(\eM,\eN)=\emptyset$, 
by Proposition~\ref{prop42}. Hence, $T^{-1}\not\in\Mod(\fD)$.
By Proposition~\ref{prop51}, $T\Alg(\fD)$ is a proper subspace of $\Mod(\fD)$.

A direct sum of copies of this example, with $2^{-j}H$ in place of $H$ in the $j$-th copy, gives the same conclusion
for $\dim\eL=\infty$. Let $T_j$ be the operator obtained in the $j$-th copy. The operators $T_j$ and $T_j^{-1}$ are uniformly bounded, 
and $\bigoplus_j2^{-j}H$ is compact.\qed
\end{example}

\subsection*{Acknowledgments}
The author thanks Professors V.~M\"{u}ller, M.~Kandi\'{c}, and B.~Kuzma for discussing the problem addressed in this paper.

\subsection*{AI Declaration}
During the preparation of this work, the author used Claude (Anthropic, model Claude Opus 5.5). The AI assistant
was used to search for and check relevant literature, to propose arguments and examples, and for language and
\LaTeX\ editing. The author verified all proofs, revised the manuscript, and takes full responsibility for its content.



\begin{thebibliography}{99}

\bibitem{At}
	F. V. Atkinson,
	\textit{The normal solubility of linear equations in normed spaces},
	Mat. Sbornik N.S. \textbf{28(70)} (1951), 3--14 (in Russian).

\bibitem{ASS}
	J. Avron, R. Seiler, B. Simon,
	\textit{The index of a pair of projections},
	J. Funct. Anal. \textbf{120} (1994), 220--237.

\bibitem{BS}
	A. B\"ottcher, I. M. Spitkovsky,
	\textit{A gentle guide to the basics of two projections theory},
	Linear Algebra Appl. \textbf{432} (2010), 1412--1459.

\bibitem{BK22}
	J. Bra\v{c}i\v{c}, M. Kandi\'c,
	\textit{On the normalizer of the reflexive cover of a unital algebra of linear transformations},
	Linear Algebra Appl. \textbf{653} (2022), 207--230.

\bibitem{BK26}
	J. Bra\v{c}i\v{c}, M. Kandi\'c,
	\textit{Similarities of subspace lattices in Banach spaces},
	J. Math. Anal. Appl. \textbf{564} (2026), 130857.

\bibitem{Con}
	J. B. Conway,
	\textit{A course in functional analysis}, 2nd ed.,
	Graduate Texts in Mathematics 96, Springer-Verlag, New York, 1990.

\bibitem{Hal}
	P. R. Halmos,
	\textit{Two subspaces},
	Trans. Amer. Math. Soc. \textbf{144} (1969), 381--389.

\bibitem{LT}
	M. Lauzon, S. Treil,
	\textit{Common complements of two subspaces of a Hilbert space},
	J. Funct. Anal. \textbf{212} (2004), 500--512.

\bibitem{LWD}
	Y. Li, Y. Wang, H. Du,
	\textit{An interpolation property of reflections involving orthogonal projections},
	Oper. Matrices \textbf{16} (2022), 265--277.

\bibitem{LP}
	W. E. Longstaff, O. Panaia,
	\textit{On complementary subspaces of Hilbert space},
	Proc. Amer. Math. Soc. \textbf{126} (1998), 3019--3026.

\bibitem{Pap}
	M. Papadakis,
	\textit{On isomorphisms between certain non-CSL algebras},
	Proc. Amer. Math. Soc. \textbf{119} (1993), 1157--1164.

\bibitem{ZH}
	L. Zhang, G. Hai,
	\textit{The properties of the product of two orthogonal projections},
	Oper. Matrices \textbf{19} (2025), 165--182.

\end{thebibliography}
\end{document}